\documentclass{article}

\usepackage{amsmath, amssymb}
\usepackage{graphicx, wrapfig} 
\usepackage{amsthm}
\usepackage{url}
\usepackage{color}
\usepackage{mathrsfs}
\usepackage{bm}

\theoremstyle{plain}
\newtheorem{theorem}{Theorem}

\theoremstyle{definition}
\newtheorem{definition}{Definition}
\newtheorem{example}{Example}

\newtheorem{proposition}{Proposition}
\newtheorem{lemma}{Lemma}
\newtheorem{corollary}{Corollary}

\newtheorem{remark}{Remark}

\usepackage{color}

\title{Crossing matrix and a polynomial invariant of braid systems up to Hurwitz equivalence}
\author{Ayaka Shimizu\thanks{Center for Soft Matter Physics, Ochanomizu University, 2-1-1, Otsuka, Bunkyo-ku, Tokyo, 112-8610, Japan. Email: shimizu.ayaka@ocha.ac.jp, shimizu1984@gmail.com} 
and Yoshiro Yaguchi\thanks{Maebashi Institute of Technology, 460-1, Kamisadori, Maebashi, Gunma, 371-0816, Japan. Email: y.yaguchi@maebashi-it.ac.jp}}
\date{\today}

\begin{document}

\maketitle

\begin{abstract}
We study the crossing matrix of a braid and introduce a polynomial invariant for braid systems that is invariant under Hurwitz equivalence. 
As an application to the study of surface braids and surface links, we also define an invariant that can be used as an indicator of the necessity of Euler fusion or fission between braid systems.
\end{abstract}



\section{Introduction}
\label{section-intro}

The Hurwitz action plays an important role in the study of surface braids and surface links. 
Let $G$ be a group, and let $G^n$ denote the $n$-fold direct product of $G$. 
For $g, h \in G$, we define an operation $*$ on $G$ by $g*h=h^{-1}gh\in G$. 
The {\it Hurwitz action} is defined by a natural right action of $B_n$ on $G^n$ as
\begin{align*}
& (g_1, g_2, \dots , g_n)\cdot \sigma_i= (g_1, \dots , g_{i-1}, g_{i+1}, g_i * g_{i+1}, g_{i+2}, \dots , g_n),\\
& (g_1, g_2, \dots , g_n)\cdot {\sigma_i}^{-1}= (g_1, \dots , g_{i-1}, g_{i+1} * (g_i^{-1}), g_i, g_{i+2}, \dots , g_n)
\end{align*}
for each $i \in \{ 1, 2, \dots , n-1 \}$, where $\sigma_i$ is the $i^{th}$ standard generator of $B_n$. Two elements $\vec{g}$ and $\vec{g'}$ of $G^n$ are said to be {\it Hurwitz equivalent} and are denoted by $\vec{g} \stackrel{{ \mathrm{Hur}}}{\backsim} \vec{g}'$ if they are related by a Hurwitz action. \medskip

When $G$ is a braid group $B_m$, a complete classification of elements of $G^n$ up to Hurwitz equivalence yields a complete classification of {\it surface braids}, that is the two-dimensional version of the classical braids, of degree $m$ with $n$ branch points up to ``equivalence" (\cite{K1996, K-b}). 
From a computational viewpoint, determining Hurwitz equivalence is difficult. 
In fact, the Hurwitz equivalence problem is undecidable in general (\cite{L}). 
This motivates the search for effective, computable invariants of Hurwitz equivalence classes. \medskip

A {\it crossing matrix}, introduced in \cite{Bu}, is a matrix whose entries represent the crossing information between the strands of a braid. 
Using the crossing matrix, several conjugacy invariants of braids were proposed in \cite{AS}. 
It is natural to extend these invariants to {\it braid systems} $\vec{b} \in (B_m)^n$. 
In Section \ref{section-const}, we construct new invariants of braid systems. 
In particular, we define the {\it characteristic polynomial} $P(\vec{b})$ of a braid system $\vec{b} \in (B_m)^n$ and prove that $P(\vec{b})$ is invariant under Hurwitz equivalence. 
The following theorem is proved in Section \ref{section-const}. 

\medskip
\begin{theorem}
Let $\vec{b}, \vec{b'}\in {(B_m)}^n$. 
If $\vec{b} \stackrel{{ \mathrm{Hur}}}{\backsim} \vec{b}'$, then $P(\vec{b})=P(\vec{b}')$. 
\label{thm-main}
\end{theorem}
\medskip

\noindent For example, let $\vec{b}=(\sigma_1 \sigma_2 \sigma_3^{-1}, \sigma_3, \sigma_2^{-1}, \sigma_1^{-1}),\ \vec{b'}=(\sigma_1 \sigma_2^{-1} \sigma_3, \sigma_3^{-1}, \sigma_2, \sigma_1^{-1})\in {(B_4)}^4$. Their characteristic polynomials are given by $P(\vec{b})=x^{16}-5x^{14}+10x^{12}-10x^{10}+5x^8-x^6$ and $P(\vec{b'})=x^{16}-9x^{14}+8x^{13}+18x^{12}-24x^{11}-10x^{10}+24x^9-3x^8-8x^7+3x^6$. 
Hence, by the contrapositive of Theorem \ref{thm-main}, we conclude that $\vec{b} \stackrel{{ \mathrm{Hur}}}{\not\backsim} \vec{b}'$. 
Note that $\vec{b}$ and $\vec{b'}$ have the same ``trace product'' and the same ``monodromy group'' (see Example \ref{ex-H} for more details). 
This example shows that $P(\vec{b})$ is an effective invariant in practice. 
Moreover, $P(\vec{b})$ is easy to compute. \medskip

The characteristic polynomial $P(\vec{b})$ of a braid system $\vec{b} \in (B_m)^n$ admits a factorization into $mn$ linear expressions as $P(\vec{b})=(x-a_1) (x-a_2) \dots (x-a_{mn})$ with real constant terms $a_1, a_2, \dots , a_{mn} \in \mathbb{R}$ (Corollary \ref{cor-factor}). 
Let $E(\vec{b})$ be the multiset of the constant terms $a_1, a_2, \dots , a_{mn}$ except for $0, \pm 1$. 
We apply this $E(\vec{b})$ to the study of (orientable and oriented) {\it surface links}, that is, a closed oriented surfaces embedded in ${\mathbb R}^4$. 
Kamada proved that every surface link can be obtained as the closure of some surface braid (\cite{K-b}). 
He also proved that if two braid systems $\vec{b}$ and $\vec{b'}$ represent ``equivalent'' surface links, then they are related by a finite sequence of a Hurwitz action, global conjugation, 
stabilization/destabilization, and Euler fusion/fission (see Sections \ref{section-const} and \ref{section-surface-link} for details). 
The following theorem is shown in Section \ref{section-const}. 

\medskip 
\begin{theorem}
Let $\vec{b}$, $\vec{b'}$ be braid systems. If they are related by a finite sequence of a Hurwitz action, global conjugation, and stabilization/destabilization, then $E(\vec{b})=E( \vec{b'})$.
\label{thm-main-2}
\end{theorem}
\medskip 

\noindent Take two braid systems $\vec{b}$ and $\vec{b'}$ which represent equivalent surface links. 
As a consequence of Theorem \ref{thm-main-2}, if $E(\vec{b}) \neq E( \vec{b'})$, then any sequence of a Hurwitz action, global conjugation, 
stabilization/destabilization, and Euler fusion/fission  between $\vec{b}$ and $\vec{b'}$ must contain an Euler fusion/fission. 
Thus, $E(\vec{b})$ can be used as an indicator for the necessity of Euler fusion or fission for a given pair of braid systems representing equivalent surface links, as shown in Example \ref{ex-E} in Section \ref{section-surface-link}. \\

The rest of the paper is organized as follows. 
In Section \ref{section-CM}, we review braids and the crossing matrix. 
In Section \ref{section-conj}, we recall the conjugacy invariants derived from the crossing matrix. 
In Section \ref{section-poly}, we explore the characteristic polynomial of braids. 
In Section \ref{section-Hur}, we review the Hurwitz action and Hurwitz equivalence. 
In Section \ref{section-const}, we construct Hurwitz invariants on braid systems and prove Theorems \ref{thm-main} and \ref{thm-main-2}. 
In Section \ref{section-braided-surface}, we discuss an application to braided surfaces and surface braids. 
In Section \ref{section-surface-link}, we recall surface links, and see an example of Theorem \ref{thm-main-2}.

\section{Braids and the crossing matrix}
\label{section-CM}

In this section, we briefly review braids, braid groups, and the crossing matrix introduced in \cite{Bu}. \medskip

Throughout this paper, we work in the PL category. 
A {\it geometric $m$-braid} ($m \geq 2$) consists of $m$ strands attached to two horizontal bars in $\mathbb{R}^3$, where each strand runs monotonically from the upper bar to the lower bar with respect to the vertical coordinate. 
Two geometric $m$-braids are said to be {\it equivalent} if they are sent to each other by an ambient isotopy of ${\mathbb{R}}^3$ which fixes the two horizontal bars pointwise. 
For geometric $m$-braids $b$ and $c$, we define the {\it braid product} $bc$ by the geometric $m$-braid obtained by placing $b$ above $c$ so that the endpoints of $b$ on the lower bar coincide with those of $c$ on the upper bar. \\

\begin{figure}[ht]
\centering
\includegraphics[width=1.5cm]{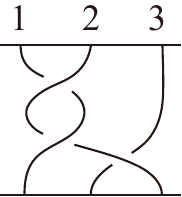}
\caption{A braid diagram. }
\label{fig-BD}
\end{figure}

An $m$-{\it braid diagram} $B$ of a geometric $m$-braid $b$ is a projection of 
$b$ onto ${\mathbb R}^2$ in general position (so that all crossings are transverse double points), equipped with over/under information at each crossing, as shown in Figure \ref{fig-BD}. 

\begin{figure}[ht]
\centering
\includegraphics[width=6cm]{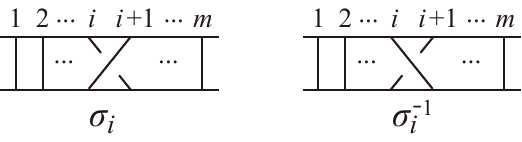}
\caption{The generators $\sigma_i$ and $\sigma_i^{-1}$. }
\label{fig-sigma}
\end{figure}

In this paper, the equivalence class of a geometric $m$-braid is simply called an $m$-{\it braid}. Let $B_m$ be the $m$-braid group, that is, the set of
$m$-braids with the group operation naturally induced by the braid product. 
For $1\leq i\leq m-1$, $\sigma_i$ is the $m$-braid which has the diagram depicted in Figure~\ref{fig-sigma} and we call it the $i^{th}$ {\it standard generator} of $B_m$. 
It is well-known that $B_m$ has the following presentation (\cite{Artin-1, Artin-2}):
\begin{equation*}
B_m = 
 \left<\sigma_1, \dots, 
\sigma_{m-1}\ \left|\ \begin{array}{rrr}\sigma_i\sigma_j\sigma_i =\sigma_j\sigma_i\sigma_j
\         (|i-j|=1),\\
\sigma_i\sigma_j=\sigma_j\sigma_i \ (|i-j|>1)
                  \end{array}\right.\right> .
\end{equation*}
For example, the braid which has the diagram in Figure \ref{fig-BD} is expressed as $\sigma_1 \sigma_1 \sigma_2^{-1}\in B_3$, or equivalently $\sigma_1^2 \sigma_2^{-1}\in B_3$. \\

\begin{figure}[ht]
\centering
\includegraphics[width=4.5cm]{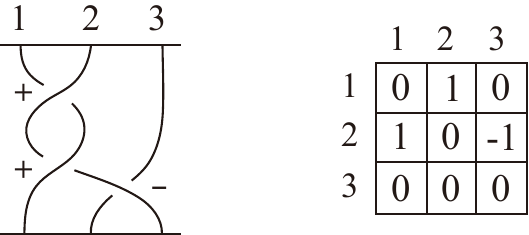}
\caption{A braid diagram and its crossing matrix. }
\label{fig-3}
\end{figure}

We recall the definition of the crossing matrix introduced in \cite{Bu}. 
For a braid or braid diagram, we call the strand that has the upper endpoint at the $k^{th}$ position from the left the {\it $k^{th}$ strand}. 
The {\it crossing matrix} was defined by Burillo, Gutierrez, Krsti\'{c} and Nitecki in \cite{Bu} (see also \cite{Gu}) as follows. 
Let $B$ be an $m$-braid diagram. 
The crossing matrix $C(B)$ is an $m \times m$ matrix with zero diagonal entries. 
Its $(i,j)$-entry is defined as the number of positive crossings minus the number of negative crossings between the $i^{th}$ and $j^{th}$ strands at which the $i^{th}$ strand passes over the $j^{th}$ strand. 
(See Figure \ref{fig-3}.) 
Since the crossing matrix is unchanged by the relation $\sigma_i\sigma_i^{-1}=\sigma_i^{-1}\sigma_i=id_{B_m}$ and two relations in the presentation of $B_m$, the crossing matrix does not depend on the choice of a braid diagram $B$ for each braid $b$. 
This enables us to define the crossing matrix $C(b)$ of a braid $b$ as $C(b)=C(B)$ for any diagram $B$ of $b$. \\

For each $m$-braid $b$, the {\it braid permutation of} $b$ is the permutation on $(1, 2, \dots , m)$ that is associated to the correspondence from the upper endpoints to the lower endpoints of the strands of $b$. 
For example, the braid permutation of the braid in Figure \ref{fig-3} is $\binom{1 \ 2 \ 3}{1 \ 3 \ 2}$ since the $1^{st}$, $2^{nd}$, $3^{rd}$ strands have the endpoints at the $1^{st}$, $3^{rd}$, $2^{nd}$ position at the bottom. 
A {\it pure braid} is a braid whose braid permutation is trivial. 
The following proposition was proved in \cite{Bu} for pure braids. 

\medskip
\begin{proposition}[\cite{Bu}]
The crossing matrix of a pure braid is symmetric. 
\label{prop-sym}
\end{proposition}
\medskip

\section{Conjugacy invariants of braids}
\label{section-conj}

Here we review and explore conjugacy invariants of braids that are proposed in \cite{AS}. 

\medskip
\begin{definition}
Two braids $b, b' \in B_m$ are said to be {\it conjugate} and denoted by $b \stackrel{{ \mathrm{conj}}}{\backsim} b'$ if $b'=b * a \ (=a^{-1}ba)$ for some $a \in B_m$.
\end{definition}
\medskip

\begin{definition}
Two $m \times m$ matrices $M$ and $N$ are said to be {\it permutation equivalent}, denoted by $M  \stackrel{{ \mathrm{perm}}}{\backsim} N$, if $N$ is obtained from $M$ by applying the same permutation to both its rows and columns. 
Equivalently, $N=P^{-1}MP$ for some $m\times m$ permutation matrix $P$. 
\end{definition}
\medskip

\begin{example}
The matrices 
$\begin{bmatrix}
1 & 2 & 3 \\
4 & 5 & 6 \\
7 & 8 & 9
\end{bmatrix}$
and 
$\begin{bmatrix}
5 & 4 & 6 \\
2 & 1 & 3 \\
8 & 7 & 9
\end{bmatrix}$
are permutation equivalent since they are related by the permutation switching the 1st and 2nd rows and columns. 
\end{example}
\medskip

\noindent Take two braids $b, b'\in B_m$. It was shown in \cite{AS} that $C(b) \stackrel{{ \mathrm{perm}}}{\backsim} C(b')$ if $b$ is {\it pure} and $b \stackrel{{ \mathrm{conj}}}{\backsim} b'$. 
In general, let $r$ be the order of the braid permutation of a braid $b$. 
Then $b^r$ is a pure braid and we obtain the following proposition. 

\medskip 
\begin{proposition}[\cite{AS}]
For two braids $b, b'\in B_m$, if $b \stackrel{{ \mathrm{conj}}}{\backsim} b'$, then $C(b^r) \stackrel{{ \mathrm{perm}}}{\backsim} C((b')^r)$, where $r$ is the order of the braid permutation of $b$. 
 
\label{prop-conj-perm}
\end{proposition}
\medskip

\begin{example}
For the two braids $b=\sigma_1 \sigma_2 \sigma_3^{-1}$ and $b'= \sigma_1 \sigma_2^{-1} \sigma_3 \in B_4$ in Figure \ref{fig-4}, the orders of their braid permutations are both 4, and we obtain the matrices 
\begin{align*}
C(b^4)=
\begin{bmatrix}
0 & 0 & 1 & 0 \\
0 & 0 & 0 & 1 \\
1 & 0 & 0 & 0 \\
0 & 1 & 0 & 0 
\end{bmatrix}, \ 
C((b')^4)=
\begin{bmatrix}
0 & 1 & -1 & 1 \\
1 & 0 & 1 & -1 \\
-1 & 1 & 0 & 1 \\
1 & -1 & 1 & 0 
\end{bmatrix}.
\end{align*}
Since $C(b^4)$ is not permutation equivalent to $C((b')^4)$, we conclude that $b$ is not conjugate to $b'$ by the contrapositive of Proposition \ref{prop-conj-perm}. 
\begin{figure}[ht]
\centering
\includegraphics[width=5.5cm]{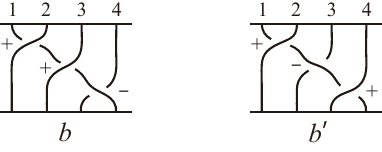}
\caption{Braids $b=\sigma_1 \sigma_2 \sigma_3^{-1}$ and $b'= \sigma_1 \sigma_2^{-1} \sigma_3$. }
\label{fig-4}
\end{figure}
\label{ex-4}
\end{example}
\medskip

\noindent In Example \ref{ex-4}, it was easy to see that $C(b^4) \stackrel{{ \mathrm{perm}}}{\not\backsim}   C((b')^4)$ since the entries are not the same. 
In some cases, it is not easy to determine whether two matrices are permutation equivalent or not. 
By the well-known fact that the rank, determinant, and characteristic polynomial are invariant under permutation equivalence for square matrices, the following proposition follows from Proposition \ref{prop-conj-perm}. 

\medskip
\begin{proposition}[\cite{AS}]
Let $r$ be the order of the braid permutation of a braid $b$. 
The rank, determinant, characteristic polynomial, and eigenvalues of $C(b^r)$ are invariant under conjugation. 
\label{prop-rank}
\end{proposition}
\medskip

\begin{example}
For the pair of conjugate braids $b=\sigma_3 \sigma_1^{-1} \sigma_4$ and $b'= \sigma_4 \sigma_3 \sigma_1^{-1}$ shown in Figure \ref{fig-5}, we have $\det (C(b^6)) = \det (C((b')^6))=-144$. 
The characteristic polynomials of $C(b^6)$ and $C((b')^6)$ are both $x^5-21x^3-16x^2+108x+144=(x-4)(x-3)(x+2)^2(x+3)$. 
\begin{figure}[ht]
\centering
\includegraphics[width=6cm]{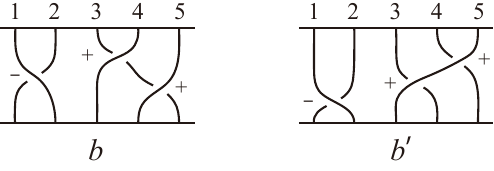}
\caption{Conjugate braids $b=\sigma_3 \sigma_1^{-1} \sigma_4$ and $b'= \sigma_4 \sigma_3 \sigma_1^{-1}$. }
\label{fig-5}
\end{figure}
\label{ex-poly5}
\end{example}
\medskip

\begin{example}
For the braids $b=\sigma_1 \sigma_2 \sigma_3^{-1}$ and $b'= \sigma_1 \sigma_2^{-1} \sigma_3$ in Example \ref{ex-4}, we obtain $\det (C(b^4))=1$, $\det (C((b')^4))=-3$. 
The characteristic polynomials of $C(b^4)$, $C((b')^4)$ are $x^4-2x^2+1=(x+1)^2 (x-1)^2$, $x^4-6x^2+8x-3={(x-1)}^3(x+3)$, respectively. 
\label{ex-poly4}
\end{example}
\medskip 

\noindent In the following proposition, further conjugacy invariants are obtained from Proposition \ref{prop-conj-perm} in the same way as in \cite{AY-d} as permutation-equivalence invariants. 

\medskip
\begin{proposition}
Let $r$ be the order of the braid permutation of a braid $b$. 
The multiset $S^R(b)$ (resp. $S^C(b)$) consisting of the multisets of the entries in each row (resp. column) of $C(b^r)$ is invariant under conjugation. 
The multiset $S(b)$ of the entries of $C(b^r)$ is also invariant under conjugation. 
\label{prop-multi}
\end{proposition}
\medskip

\begin{example}
For the braid $b$ shown in Figure \ref{fig-3}, we have $S^R(b)= S^C(b)= \{ \{ 0,1,1 \}, \{ 1,0,-1 \}, \{ 1, -1, 0 \} \}$ and $S(b)= \{ -1, -1, 0, 0, 0, 1, 1, 1, 1 \}$.
\end{example}
\medskip

\noindent From the next section, we focus on the characteristic polynomial of $C(b^r)$, which is reasonable to deal with and contains the information of the determinant and eigenvalues. 

\section{The characteristic polynomial of braids}
\label{section-poly}

\noindent We investigate the characteristic polynomial of the crossing matrix $C(b^r)$ for a braid $b$ with order of the braid permutation $r$. 

\medskip
\begin{definition}
The {\it characteristic polynomial} of $b \in B_m$, denoted by $P(b)$, is the characteristic polynomial of the crossing matrix $C(b^r)$, where $r$ is the order of the braid permutation of $b$. 
That is, $P(b)= \det \left( xI - C(b^r) \right)$, where $I$ is the $m \times m$ identity matrix. 
\end{definition}
\medskip 

\noindent As shown in Proposition \ref{prop-rank}, $P(b)$ is invariant under conjugation. 
The following propositions hold (see Examples \ref{ex-poly5} and \ref{ex-poly4}). 

\medskip 
\begin{proposition}
For each $b \in B_m$, $P(b)$ is in the following form: 
\begin{align*}
P(b)= x^m + c_{m-2}x^{m-2} + c_{m-3}x^{m-3} + \dots + c_1 x +c_0,
\end{align*}
where $c_{m-2}, c_{m-3}, \dots , c_1 , c_0 \in \mathbb{Z}$. 
\end{proposition}
\medskip 

\begin{proof}
The coefficient of $x^m$ in $P(b)$ is one by the property of the characteristic polynomial. 
Let $r$ be the order of the braid permutation of $b$. 
The coefficient of $x^{m-1}$ is equal to the trace of $-C(b^r)$, which is zero because the crossing matrix has zero diagonal. 
Since $C(b^r)$ is an integer matrix, $c_{m-2}, c_{m-3}, \dots , c_1 , c_0 \in \mathbb{Z}$. 
\end{proof}

\medskip
\begin{proposition}
For each $b \in B_m$, $P(b)$ factors as $P(b)=(x-a_1)(x-a_2) \dots (x-a_m)$ with real numbers $a_1, a_2, \dots , a_m \in \mathbb{R}$ such that $a_1 +a_2 + \dots + a_m =0$. 
\label{prop-b-factor}
\end{proposition}
\medskip

\begin{proof}
Let $r$ be the order of the braid permutation of $b$. 
As observed in \cite{AS}, the $m \times m$ matrix $C(b^r)$ has real eigenvalues, say $a_1, a_2, \dots , a_m \in \mathbb{R}$ since $b^r$ is pure and $C(b^r)$ is symmetric by Proposition \ref{prop-sym}. 
Note that $a_1 + a_2 + \dots + a_m$ is equal to the trace of $C(b^r)$, which is zero. 
\end{proof}
\medskip

\noindent In the rest of this section, we observe some formulae of the characteristic polynomials of braids. 

\medskip
\begin{example}
For each $i\in \{1,2,\dots, m-1\}$, we have $P(\sigma_i)=P(\sigma_i^{-1})=x^{m-2}(x+1)(x-1)$, where $\sigma_i$ is the $i^{th}$ standard generator of $B_m$.
\label{ex-sigma}
\end{example}

\begin{proof}The order of the braid permutation of ${\sigma_i}^{\varepsilon}$ is two for each $\varepsilon\in \{\pm 1\}$, and $C(\sigma_i^2)$ (resp. $C(\sigma_i^{-2})$) is an $m \times m$ matrix such that the $(i, i+1)$- and $(i+1, i)$-entries are 1 (resp. $-1$) and the others are zero. 
Thus, we obtain $P(\sigma_i)=P(\sigma_i^{-1})=x^{m-2}(x+1)(x-1)$ for any  $\sigma_i, \sigma_i^{-1} \in B_m$. 
\end{proof}

\begin{figure}[ht]
\centering
\includegraphics[width=2cm]{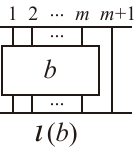}
\caption{The braid $\iota(b)$. }
\label{fig-plus}
\end{figure}

\medskip 
\begin{lemma}
Let $b \in B_m$. 
Let $\iota(b)$ be the $(m+1)$-braid that has the same word as $b$ (see Figure \ref{fig-plus}). 
Then $P(\iota(b))=xP(b)$. 
\label{lem-plus}
\end{lemma}

\medskip 
\begin{proof}
Let $r$ be the order of the braid permutation of $b$. 
Then $\iota(b)$ has the same order $r$ of the braid permutation, and we have $(\iota(b))^r=\iota(b^r)$. 
The $(m+1)\times (m+1)$ crossing matrix $C(\iota(b^r))$ coincides with the $m \times m$ crossing matrix $C(b^r)$ in the first $m$ rows and columns, and 
all the entries in the $(m+1)^{th}$ row and column are 0. 
Hence, $P(\iota(b))= \det \left( x I_{m+1} - C ((\iota(b))^r) \right) = \det \left( x I_{m+1} - C (\iota(b^r)) \right) = x \det \left( x I_{m} - C (b^r) \right) = x P(b)$, where $I_k$ denotes the $k \times k$ identity matrix. 
\end{proof}
\medskip

\noindent Let $id\in B_m$ be the trivial braid. We have $P(id)=x^m$. By considering ``weaving braids'', we obtain the following proposition. 

\medskip 
\begin{proposition} (1) When $m=2$, there are not nontrivial $m$-braids $b\in B_m$ such that $P(b)=x^m$ and (2) for each $m\geq 3$, there are nontrivial $m$-braids $b$ such that $P(b)=x^m$. 
\end{proposition}
\medskip

\begin{proof}(1) Let $b= \sigma_1^l \in B_2$, where $l$ is an integer. 
Since the order $r$ of the braid permutation of $b$ is 2 (resp. 1) when $l$ is odd (resp. even), $b^r$ is in the form of $b^r=(\sigma_1)^{2k}$, where $k=l$ (resp. $k=\frac{l}{2}$) when $l$ is odd (resp. even). 
The $(1,2)$- and $(2,1)$-entries of $C(b^r)$ are both $k$, and $P(b)=x^2-k^2$. 
We have $P(b)=x^2$ if and only if $k=0$, that is, $l=0$. \\
(2) Let $n\geq 3$ be a positive odd number. 
Let $b= \sigma_1 \sigma_2^{-1} \sigma_3 \dots \sigma_{n-1}^{-1}$, namely a weaving braid $W(n,1)$. 
Then $b$ has the order of the braid permutation as $n$. 
As shown in \cite{ASA}, each pair of strand of $b^n (=W(n,n))$ has exactly two crossings, one positive crossing and one negative crossing, where one strand is over the other strand at both crossings\footnote{When $n$ is a positive even number with $n\geq 4$, one of the strands has one over-crossing and one under-crossing for each pair of strands of $W(n,n)$. See Section 2.1 of \cite{ASA} for more details.}. 
Hence, the crossing matrix of $b^n$ is the $n \times n$ zero matrix $O$. 
Therefore, $P(b)=P(b^n)= \det ( xI-O)=x^n$. 
We also have $P(\iota (b))=P(\iota(b^n)) =xP(b) =x^{n+1}$ by Lemma \ref{lem-plus}. 
\end{proof}
\medskip 

\noindent We provide formulae of the characteristic polynomial for some positive pure braids\footnote{The characterization of the crossing matrix of a positive pure braid is an open problem. See \cite{Bu}, \cite{Gu}, \cite{YAY}, \cite{AY-5}.}. 

\medskip 
\begin{example}
For $m>2$, let $b_m= \sigma_1 \sigma_2 \dots \sigma_{m-2} \sigma_{m-1}^2 \sigma_{m-2} \dots \sigma_2 \sigma_1 \in B_m$ (see Figure \ref{fig-p5}). 
Then $P(b_m)= x^m -(m-1) x^{m-2}$. 
\begin{figure}[ht]
\centering
\includegraphics[width=1.8cm]{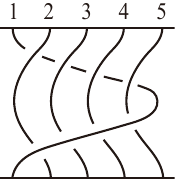}
\caption{The braid $b_5$ has the characteristic polynomial $P(b_5)=x^5-4x^3$. }
\label{fig-p5}
\end{figure}
\label{ex-bm}
\end{example}
\medskip 

\begin{proof}
Since $b_m$ is a pure braid, we have $P(b_m) = \det \left( xI-C(b_m) \right)$. 
By the Laplace expansion along the $m^{th}$ column and then the $(m-1)^{th}$ row, we have
\begin{align*}
& \begin{vmatrix}
x & -1 & -1 & \dots & -1 \\
-1 & x & 0 & \dots & 0 \\
-1 & 0 & x & \dots & 0\\
 \vdots & \vdots & \vdots & \ddots & \vdots \\
 -1 & 0 & 0 & \dots & x
\end{vmatrix}
= (-1)^m 
\begin{vmatrix}
-1 & x & 0 & \dots & 0 \\
-1 & 0 & x & \dots & 0 \\
 \vdots &  \vdots & \vdots & \ddots & \vdots \\
 -1 & 0 & 0 & \dots & x \\
 -1 & 0 & 0 & \dots & 0
\end{vmatrix}
+x 
\begin{vmatrix}
x & -1 &  \dots & -1 \\
-1 & x &  \dots & 0 \\
 \vdots & \vdots &  \ddots & \vdots \\
 -1 & 0 &  \dots & x
\end{vmatrix} \\
& =(-1)^m (-1)^{m-1}
\begin{vmatrix}
x & 0 &  \dots & 0 \\
0 & x &  \dots & 0 \\
 \vdots & \vdots &  \ddots & \vdots \\
 0 & 0 &  \dots & x
\end{vmatrix} 
+ xP(b_{m-1}) =-x^{m-2} +x P(b_{m-1}).
\end{align*}
By induction, the formula $P(b_m)=x^m-(m-1)x^{m-2}$ is proved. 
\end{proof}
\medskip 

\noindent For positive pure 3-braids, the following proposition holds. 

\medskip 
\begin{proposition}
Let $b$ be a positive pure 3-braid. 
The characteristic polynomial is in the form $P(b)=x^3+a_1 x +a_0$, where $a_1$ is a non-positive integer and $a_0$ is an even non-positive integer. 
Moreover, $a_0$ is non-zero if and only if every pair of strands has a crossing in any diagram of $b$. 
\end{proposition}
\medskip 

\begin{proof}
Let $2k, 2l, 2j$ be the number of crossings between the $1^{st}$ and $2^{nd}$, $1^{st}$ and $3^{rd}$, $2^{nd}$ and $3^{rd}$ strands of a diagram of $b$. 
Since $C(b)$ is a symmetric matrix,  
\begin{align*}
P(b)= \det (xI-C(b)) = 
\begin{vmatrix}
x & -k & -l \\
-k & x & -j \\
-l & -j & x
\end{vmatrix}
=x^3 -(k^2+l^2+j^2)x -2klj.
\end{align*}
The constant term is non-zero if and only if all of $k, l, j$ are non-zero. 
\end{proof}
\medskip

\section{The Hurwitz action and its properties}
\label{section-Hur}

In this section, we briefly review the Hurwitz action and the Hurwitz equivalence. 
Let $G$ be a group and let $G^n$ be the $n$-fold direct product of $G$. 
We restate the definition of the Hurwitz action and the Hurwitz equivalence for convenience. 

\medskip
\begin{definition}The Hurwitz action of $B_n$ on $G^n$ is the right action defined by
\begin{align*}
(g_1, g_2, \dots , g_n)\cdot \sigma_i=(g_1, \dots , g_{i-1}, g_{i+1}, g_i * g_{i+1}, g_{i+2}, \dots , g_n)\ {\rm and}\\
(g_1, g_2, \dots , g_n)\cdot {\sigma_i}^{-1}=(g_1, \dots , g_{i-1}, g_{i+1} * (g_i^{-1}), g_i, g_{i+2}, \dots , g_n)
\end{align*}
for each $i \in \{ 1, 2, \dots , n-1 \}$. 
Two elements $\vec{g}=(g_1, \dots , g_n)$ and $\vec{h}=(h_1, \dots , h_n)$ of $G^n$ are {\it Hurwitz equivalent} and denoted by $\vec{g} \stackrel{{ \mathrm{Hur}}}{\backsim} \vec{h}$ if they are sent to each other by the Hurwitz action of $B_n$, i.e., there exists an element $\beta\in B_n$ such that $\vec{g}\cdot \beta=\vec{h}$.
\end{definition}
\medskip

\begin{definition}
For an element $\vec{g}=(g_1, \dots , g_n)\in G^n$,
\begin{itemize}
\item[(1)] the {\it trace product of $\vec{g}$} is the 
product $g_1\cdots g_n$ in $G$ and denoted by $tr(\vec{g})$, and
\item[(2)] the {\it monodromy group of $\vec{g}$} is the subgroup of $G$ generated by $g_1,\dots ,g_n$ and denoted by $\langle \vec{g}\rangle$.
\end{itemize}
\medskip
\end{definition}

\noindent The following properties are obtained from the definition of the Hurwitz action (see, e.g., \cite{EB}, \cite{L}). 

\medskip
\begin{proposition}
Let $\vec{g}=(g_1, \dots , g_n)$, $\vec{h}=(h_1, \dots , h_n) \in G^n$. 
If $\vec{g} \stackrel{{ \mathrm{Hur}}}{\backsim} \vec{h}$, then 
\begin{itemize}
\item[(A)] $tr(\vec{g})=tr(\vec{h})$, 
\item[(B)] $\langle \vec{g}\rangle =\langle \vec{h}\rangle$,
\item[(C)] for any $i\in \{1,\dots ,n\}$, $g_i$ and $h_{\pi(i)}$ are conjugate in $\langle \vec{g}\rangle$, where $\pi$ is the braid permutation of $\beta\in B_n$ such that $\vec{g}\cdot \beta=\vec{h}$, i.e., the number of times each conjugacy class with respect to $\langle \vec{g}\rangle$ appears in $\vec{g}$ is the same as in $\vec{h}$, and 
\item[(D)] for a group $G'$ and for a group homomorphism $f: G \to G'$, we have $f^n(\vec{g}) \stackrel{{ \mathrm{Hur}}}{\backsim} f^n(\vec{h})$, where $f^n: G^n\to {G'}^n$ is defined by $f^n(g_1,\dots ,g_n)=(f(g_1),\dots ,f(g_n))$. 
\end{itemize}
\label{prop-trace}
\label{prop-NC}
\end{proposition}
\medskip

Take a subset $X$ of $G$ which is {\it closed under conjugation in} $G$, i.e., for any
$x \in X$ and $y\in G$, $x * y \in X$. Then, the Hurwitz action of $B_n$ on $G^n$ can be
restricted to that on $X^n$. Fix an element $g\in G$. Let $P^n_g (X)$ be the set of elements
$\vec{g}\in X^n$ such that $tr(\vec{g}) = g$ and $g_i \ne id_G$ for $i \in \{1,\dots ,n\}$. When $P^n_g (X)\ne \emptyset$, by Proposition \ref{prop-NC}, the Hurwitz action of $B_n$ on $G^n$ can be restricted to that on $P^n_g (X)$.\\

In this paper, we discuss the Hurwitz equivalence for the case of $G=B_m$, namely, $G^n=(B_m)^n$. 

\section{Hurwitz invariants of braid systems}
\label{section-const}

We construct invariants of the Hurwitz equivalence classes of braid systems based on the necessary condition (C) in Proposition \ref{prop-NC}. \\

A {\it braid system of degree $m$ and length $n$} is an element $\vec{b}=(b_1, b_2, \dots , b_n)\in {(B_m)}^n$. 
Under the Hurwitz action on ${(B_m)}^n$, the conjugacy classes of the components of each braid system in ${(B_m)}^n$ are preserved up to permutation. 
Therefore, the multiset, product, or sum of the conjugacy invariants given in Propositions \ref{prop-rank} and \ref{prop-multi} is invariant under Hurwitz equivalence. 
We have the following corollary by Proposition \ref{prop-rank}. 

\medskip 
\begin{corollary}
Let $\vec{b}=(b_1, b_2, \dots , b_n) \in (B_m)^n$. 
Let $r_i$ be the order of the braid permutation of $b_i$ ($i=1, 2, \dots , n$). 
Any of the sum, product, or multiset of the ranks, determinants, and characteristic polynomials of $C(b_1^{r_1}), C(b_2^{r_2}), \dots , C(b_n^{r_n})$ is invariant under the Hurwitz action. 
The multiset of eigenvalues of \\
$C(b_1^{r_1}), C(b_2^{r_2}), \dots , C(b_n^{r_n})$ is also invariant under the Hurwitz action. 
\label{cor-H-rank}
\end{corollary}
\medskip 

\noindent The following example demonstrates a case that the multiset of characteristic polynomials of a braid system provides a more effective Hurwitz invariant than the classical necessary conditions such as the trace product and the monodromy group. 

\medskip 
\begin{example}
Let $\vec{b}=(b_1, b_2, b_3, b_4)$ and $\vec{b'}=(b'_1, b'_2, b'_3, b'_4)$ be braid systems in $(B_4)^4$ defined as follows: 
\begin{align*}
\vec{b}= & (b_1, b_2, b_3, b_4)=(\sigma_1 \sigma_2 \sigma_3^{-1}, \sigma_3, \sigma_2^{-1}, \sigma_1^{-1}), \\
\vec{b'}= & (b'_1, b'_2, b'_3, b'_4)=(\sigma_1 \sigma_2^{-1} \sigma_3, \sigma_3^{-1}, \sigma_2, \sigma_1^{-1}).
\end{align*}
The multisets of characteristic polynomials are $\{ x^4-2x^2+1, x^4-x^2, x^4-x^2, x^4-x^2 \}$ and $\{ x^4-6x^2+8x-3, x^4-x^2, x^4-x^2, x^4-x^2 \}$ by Examples \ref{ex-4} and \ref{ex-sigma}. 
Since these multisets do not coincide, we conclude from Corollary \ref{cor-H-rank} that $\vec{b} \stackrel{{ \mathrm{Hur}}}{\not\backsim} \vec{b}'$. 
It is important to note that the Hurwitz equivalence of $\vec{b}$ and $\vec{b'}$ cannot be distinguished using the necessary conditions (A) and (B) in Proposition \ref{prop-NC}; 
\begin{itemize}
\item[(A): ] $tr(\vec{b}) = tr(\vec{b'}) \ (=id)$, 
\item[(B): ] $\langle \vec{b}\rangle = \langle \vec{b'}\rangle \ (=B_4)$.
\end{itemize}
\noindent As for (D), moreover, the Hurwitz equivalence of $\vec{b}$ and $\vec{b'}$ cannot be distinguished for the following typical settings of $f$. 
\begin{itemize}
\item[(D-1): ] Set $f= \pi : B_4 \to S_4$, where $\pi$ sends each $b \in B_4$ to its braid permutation. 
Then $f^4(\vec{b})= \left( \binom{1 \ 2 \ 3 \ 4}{4 \ 1 \ 2 \ 3}, \binom{1 \ 2 \ 3 \ 4}{1 \ 2 \ 4 \ 3}, \binom{1 \ 2 \ 3 \ 4}{1 \ 3 \ 2 \ 4}, \binom{1 \ 2 \ 3 \ 4}{2 \ 1 \ 3 \ 4} \right) = f^4(\vec{b'})$, and $f^4(\vec{b}) \stackrel{{ \mathrm{Hur}}}{\backsim} f^4(\vec{b}')$.
\item[(D-2): ] Set $f: B_4 \to \mathbb{Z}$ as a map that sends each braid $b= \sigma_{i_1}^{\varepsilon_1} \sigma_{i_2}^{\varepsilon_2} \dots \sigma_{i_k}^{\varepsilon_k}$ ($i_1 , i_2 , \dots , i_k \in \{ 1, 2, 3 \}$, $\varepsilon_1, \varepsilon_2 , \dots , \varepsilon_k \in \{ \pm 1 \}$) to the sum of the exponents $\varepsilon_1 + \varepsilon_2 + \dots + \varepsilon_k$. 
Then $f^4 (\vec{b})=(1, 1, -1, -1)$, $f^4(\vec{b'})=(1, -1, 1, -1)$. 
Since $f^4 (\vec{b}) \cdot \sigma_2 = f^4 (\vec{b'})$, we obtain that $f^4(\vec{b}) \stackrel{{ \mathrm{Hur}}}{\backsim} f^4(\vec{b}')$. 
\end{itemize}
\label{ex-H}
\end{example}
\medskip

\noindent We are now ready to define the characteristic polynomial for braid systems.

\medskip
\begin{definition}
Let $\vec{b}=(b_1, b_2, \dots , b_n) \in (B_m)^n$ be a braid system. 
The {\it characteristic polynomial of $\vec{b}$}, denoted by $P(\vec{b})$, is a polynomial of degree $mn$ defined as $P(\vec{b})=P(b_1) P(b_2) \dots P(b_n)$.
\end{definition}
\medskip

\noindent We prove Theorem \ref{thm-main}, which states that if $\vec{b} \stackrel{{ \mathrm{Hur}}}{\backsim} \vec{b}'$, then $P(\vec{b})=P(\vec{b}')$.  \\

\noindent {\it Proof of Theorem \ref{thm-main}.} 
This is an immediate consequence of Corollary \ref{cor-H-rank} by taking the product of characteristic polynomials. 
\hfill$\square$  \\
\medskip 

\noindent The characteristic polynomial of a braid system has the following property. 

\medskip 
\begin{corollary}
For each $\vec{b}=(b_1, b_2, \dots , b_n) \in (B_m)^n$, the characteristic polynomial factors as 
$P(\vec{b})=(x-a_1)(x-a_2) \dots (x-a_{mn})$, where $a_1, a_2, \dots , a_{mn} \in \mathbb{R}$ and $a_1 + a_2 + \dots + a_{mn}=0$.
\label{cor-factor}
\end{corollary}
\medskip

\begin{proof}
By Proposition \ref{prop-b-factor}, each characteristic polynomial $P(b_k)$ of $b_k$ admits a factorization of the form $(x-\alpha_{k,1})(x-\alpha_{k,2})\dots (x-\alpha_{k,m})$, where $\alpha_{k,1}, \alpha_{k,2}, \dots , \alpha_{k,m} \in \mathbb{R}$ and $\alpha_{k,1} + \alpha_{k,2} + \dots + \alpha_{k,m} =0$. 
Multiplying all these characteristic polynomials gives 
\begin{align*}
P(\vec{b})=(x-\alpha_{1,1})\dots (x-\alpha_{1,m})(x-\alpha_{2,1})\dots (x-\alpha_{2,m})(x-\alpha_{n,1})\dots (x-\alpha_{n,m}).
\end{align*}
\end{proof}
\medskip

\noindent The following corollary follows from Proposition \ref{prop-multi}. 

\medskip 
\begin{corollary}
Let $\vec{b}=(b_1, b_2, \dots , b_n)\in {(B_m)^n}$. 
The multiset of the components of $S^R(b_1), S^R(b_2), \dots , S^R(b_n)$ or $S^C(b_1), S^C(b_2), \dots , S^C(b_n)$ is invariant under the Hurwitz action. 
The multiset of the components of $S(b_1), S(b_2), \dots , S(b_n)$ is also invariant under the Hurwitz action. 
\end{corollary}
\medskip

\noindent As we will see in section \ref{section-surface-link}, these Hurwitz invariants not only distinguish non–Hurwitz equivalent braid systems, but also provide an invariant regarding the necessity of ``Euler fusion or fission'' between the braid systems representing surface braids and surface links. 
To describe Markov’s theorem in dimension four in Section \ref{section-surface-link}, we prepare the following definitions and a lemma.

\medskip 
\begin{definition}
The {\it global conjugation of} $\vec{b}=(b_1, b_2, \dots ,b_n)\in {(B_m)}^n$ {\it by} $a\in B_m$ is defined by $\vec{b} * a:=(b_1 * a, \ b_2 * a, \dots ,b_n * a)\in {(B_m)}^n$. Two braid systems $\vec{b}$ and $\vec{b'}$ are {\it global conjugate} if there exists a braid $a\in B_m$ such that $\vec{b'}=\vec{b} * a$. 
\label{globalconj}
\end{definition} 
\medskip 

\noindent For each $m\in {\mathbb N}$, let $\iota : B_m \to B_{m+1}$ with $\iota(b)=b$, as defined in Lemma \ref{lem-plus}. 

\medskip 
\begin{definition}
The {\it stabilization of} $\vec{b}=(b_1, b_2, \cdots ,b_n)\in {(B_m)}^n$ is the braid system defined by $(\iota(b_1), \iota(b_2), \dots , \iota(b_n), \sigma_m, \sigma_m^{-1})\in {(B_{m+1})}^{n+2}$. 
The inverse is called a {\it destabilization} (if it exists). 
\label{stabilization}
\end{definition}
\medskip 

\noindent Observe that applying a stabilization on a braid system produces new eigenvalues $0$ or $\pm 1$ for each component by Lemma \ref{lem-plus} and Example \ref{ex-sigma}. 
Now we define an invariant for braid systems that is unchanged by any of the Hurwitz action, global conjugation, and stabilization. 

\medskip 
\begin{definition}
Let $\vec{b}=(b_1, b_2, \dots , b_n)\in {(B_m)}^n$ and let $r_i$ be the order of the braid permutation of $b_i$ for each $i\in \{1, 2, \dots ,n\}$. 
Let $E(\vec{b})$ be the multiset of eigenvalues of $C(b_1^{r_1}), C(b_2^{r_2}), \dots , C(b_n^{r_n})$ except for $0, \pm 1$. 
We refer to $E(\vec{b})$ as the {\it essential eigenvalue set of} $\vec{b}$. 
\end{definition}
\medskip

\noindent For each $\vec{b}\in {(B_m)}^n$, we can calculate $E(\vec{b})$ from the characteristic polynomial $P(\vec{b})=(x-a_1)(x-a_2) \dots (x-a_{mn})$ by taking the multiset of $a_1, a_2, \dots , a_{mn} \in \mathbb{R}$ and then eliminating $0, \pm 1$. 
\medskip 

\begin{lemma}
For a braid system $\vec{b}$, the set $E(\vec{b})$ is invariant under the transformations (I): a Hurwitz action, (I\hspace{-0.5pt}I): a global conjugation, and  (I\hspace{-0.5pt}I\hspace{-0.5pt}I): a stabilization/destabilization. 
\label{lem-123}
\end{lemma}
\medskip

\begin{proof}
\begin{itemize}
\item[(I): ] This follows from Corollary \ref{cor-H-rank}. 
\item[(I\hspace{-0.5pt}I): ] This follows from Proposition \ref{prop-rank}. 
\item[(I\hspace{-0.5pt}I\hspace{-0.5pt}I): ] For $\vec{b}=(b_1, b_2, \dots , b_n) \in (B_m)^n$ and $( \iota(b_1), \iota(b_2), \dots , \iota(b_n), \sigma_m, \sigma_m^{-1} ) \in (B_{m+1})^{n+2}$, we have $P(\iota (b_i))=x P(b_i)$ when $i=1, 2, \dots , n$ by Lemma \ref{lem-plus}. 
We have $P(\sigma_m)=P(\sigma_m^{-1})=x^{m-1}(x+1)(x-1)$ by Example \ref{ex-sigma}. 
Hence, the eigenvalues are unchanged except for $0, \pm 1$. 
\end{itemize}
\end{proof}

\section{Application to braided surfaces}
\label{section-braided-surface}

In this section, we observe that the Hurwitz invariants introduced in Section \ref{section-const} have applications to braided surfaces. We recall a braided surface (see \cite{K1996, K-b} for details). Let $D_1$ and $D_2$ be $2$-disks, and let $p_i : D_1 \times D_2 \to D_i$ be the natural projection that sends $(a_1, a_2) \in D_1 \times D_2$ to $a_i \in D_i$ for $i\in \{1,2\}$. Fix a set 
$Q_m$ of $m$ points in ${\rm Int}(D_1)$ and fix a point $y_0$ in $\partial D_2$.\medskip

A {\it braided surface} $S$ {\it of degree} $m$ is an (orientable and) oriented surface properly and locally flatly embedded in $D_1 \times D_2$ such that 
\begin{itemize}
\item[(1) ]
$p_2|_S: S \to D_2$ is a branched covering of degree $m$, i.e., there is a finite set $\Sigma_S \subset {\rm Int}(D_2)$ such that $p_2|_{S\setminus {p_2}^{-1}(\Sigma_S)}: S\setminus {p_2}^{-1}(\Sigma_S)\to D_2\setminus \Sigma_S$ is a covering map of degree $m$ and for each $x\in p_2^{-1}(\Sigma_S)\cap S$, the map $p_2|_S$ about $x$ is locally equivalent to the map $z\mapsto z^q\ (z\in {\mathbb C})$ about $0\in {\mathbb C}$ for some $q\in {\mathbb N}$,
\item[(2) ]
$\partial S\subset {\rm Int}(D_1^2)\times \partial D_2^2$ and
\item[(3) ]
$S\cap p_2^{-1}(y_0)=Q_m\times \{y_0\}$.
\end{itemize}

\noindent We call a point of the set $\Sigma_S$ in (1) a {\it branch point} of $S$. The integer $q$ in (1) is called the {\it branching index} of $x$ and is denoted by ${\rm deg}(S;x)$. For a branch point $y\in \Sigma_S$, suppose that $p^{-1}(y)\cap S$ consists of $c$ distinct points $x_1,\dots ,x_c$. It is obvious that ${\rm deg}(S;x_1)+\cdots +{\rm deg}(S;x_c)=m$. Changing indices, we may assume that ${\rm deg}(S;x_1)\geq \dots \geq {\rm deg}(S;x_c)$. The {\it branch type} of $y\in \Sigma_S$ is the $c$-tuple $({\rm deg}(S;x_1),\dots , {\rm deg}(S;x_c))$ and is denoted by ${\rm type}(S;y)$. Suppose that $|\Sigma_S|=n$ and set $\Sigma_S=\{y_1,\dots ,y_n\}$. The {\it branch type} of $S$ is the multiset $\{{\rm type}(S;y_1),\dots ,{\rm type}(S;y_n)\}$ of branch types of points in $\Sigma_S$ and is denoted by ${\rm type}(S)$.\medskip

To capture information about a braided surface beyond its branch type, we consider the following notion of a braid system. Let $S$ be a braided surface of degree $m$ with $n$ branch points. Set $\Sigma_S=\{y_1,\dots ,y_n\}$. For each loop $\gamma:(I, \partial I)\to (D_2\setminus \Sigma_S, y_0)$, we define the loop $\widetilde{\gamma}:(I, \partial I)\to (C_m, Q_m)$ by $\widetilde{\gamma}(t) = p_1(S \cap p_2^{-1} (\gamma(t)))$, where $C_m$ is the configuration space of unordered $m$ points of ${\rm Int}(D_
1)$. The fundamental group $\pi_1(C_m, Q_m)$  is isomorphic to $B_m$. Identifying $\pi_1(C_
m, Q_m)$ with $B_m$, we have a homomorphism $\rho_S: \pi_1(D_2
\setminus \Sigma_S,y_0) \to B_m$ and call it the {\it braid monodromy} of $S$. We use the same symbol $\gamma$ for the homotopy class $[\gamma]\in \pi_1(D_2\setminus\Sigma_S,y_0)$. 
A {\it Hurwitz generating system} of $\pi_1(D_2\setminus \Sigma_S, y_0)$ is an $n$-tuple $(\gamma_1,\dots , \gamma_n)$ of (the homotopy classes of) loops in $D_2\setminus \Sigma_S$ with base point $y_0$ such that each $\gamma_i$ surrounds the point $y_i\in\Sigma_S$ in a positive direction for $i\in\{1,\dots ,n\}$ and the product $\gamma_1\cdots\gamma_n$ is homotopic to the positively oriented boundary loop $\partial D_2$ based at $y_0$. Thus, we have an element $(\rho_S(\gamma_1),\dots ,\rho_S(\gamma_n))\in {(B_m)}^n$ and call it a {\it braid system} of $S$. It is known that the Hurwitz equivalence class of braid systems obtained from a braided surface is uniquely determined independently of the choice of a Hurwitz generating system. Let $f:B_m\to S_m$ be a homomorphism that sends $b\in B_m$ to its braid permutation. Note that ${\rm type}(S;y_i)$ coincides with the cycle type of $f(\rho_S(\gamma_i))$ for each $i\in \{1,\dots, n\}$, and ${\rm type}(S)$ coincides with the multiset of cycle types of $f(\rho_S(\gamma_1)),\ldots,f(\rho_S(\gamma_n))$ (see also \cite{Eis}).\medskip

To characterize the braid systems arising from braided surfaces, we
introduce the following subset $A_m$ of $B_m$. Let $D$ and $S^1$ be a 2-disk and a 1-sphere, respectively. 
Let $p: D\times S^1\to S^1$ be the natural projection. 
A {\it closed geometric $m$-braid} $\ell$ is a closed 1-manifold embedded in the solid torus $D\times S^1$ such that $p|_{\ell} : \ell\to S^1$ is a covering map of degree $m$. A closed geometric $m$-braid $\ell$ is said to be {\it completely split} if there exist $c$ mutually disjoint convex disks $N_1,\dots, N_c$ in $D$ such that every solid torus $N_i\times S^1$ contains a component of $\ell$, where $c$ is the number of connected components of $\ell$. 
For a geometric $m$-braid $b$, the {\it closure} $\widehat{b}$ is a closed braid obtained from $b$ by connecting each endpoint on the upper bar to the same position on the lower bar by disjoint arcs outside $b$ so that no new crossings are produced. 
An element $b\in B_m$ is said to be {\it completely splittable} if there is a geometric $m$-braid $b$ as a representative of it such that the closure $\widehat{b}$ is completely split. 
Let $A_m$ be the set of completely splittable non-identity braids $b\in B_m$ such that $\widehat{b}$ is a trivial link in ${\mathbb R}^3$. For example, $\sigma_i, \sigma_i^{-1}\in A_m$ for $1\leq i\leq m-1$. 
Note that $A_m$ is closed under conjugation in $B_m$.\medskip

For a braided surface $S$ of degree $m$, $\partial S$ is a closed geometric $m$-braid in the solid torus $D_1\times \partial D_2$. The {\it boundary braid of} $S$ is the $m$-braid obtained from $\partial S$ by cutting off along $p_2^{-1}(y_0) = D_1\times \{y_0\}$. Regarding $\partial D_2$ as the positively oriented boundary loop based at $y_0$, the boundary braid of $S$ coincides with $\rho_S(\partial D_2)$, that is equal to the product $\rho_S(\gamma_1)\cdots \rho_S(\gamma_n)$, where $n$ is the number of branch points of $S$ and $(\gamma_1,\dots ,\gamma_n)$ is a Hurwitz generating system of $S$.\medskip

A {\it surface braid of degree} $m$ is a braided surface $S$ of degree $m$ with $\partial S=Q_m\times \partial D_2$. Note that the boundary braid of a surface braid $S$ of degree $m$ is the trivial braid $id=id_{B_m}$.\medskip

We have the following theorem.\medskip

\begin{theorem} [\cite{K1996, K-b}]
An element $\vec{b}\in {(B_m)}^n$ is a braid system of a braided surface of degree $m$ with $n$ branch points whose boundary braid is $\delta\in B_m$ if and only if $\vec{b}\in P^n_{\delta}(A_m)$, i.e., $\vec{b}\in {(A_m)}^n$ and $tr(\vec{b})=\delta$. In particular, an element $\vec{b}\in {(B_m)}^n$ is a braid system of a surface braid of degree $m$ with $n$ branch points if and only if $\vec{b}\in P^n_{id}(A_m)$. 
\label{characterization of a braid system}
\end{theorem}
\medskip 

Two braided surfaces of degree $m$ are {\it equivalent} if they are ambient isotopic by an isotopy $\{h_u\}_{u\in [0,1]}$ of $D_1\times D_2$ such that for each $u\in [0,1]$, $h_u$ is fiber-preserving, that is, $p_2\circ h_u = \underline{h_u}\circ p_2$ for some homeomorphism $\underline{h_u} : D_2 \to D_2$ and $h_u|_{p_2^{-1}(y_0)}=id$. Take equivalent braided surfaces $S$ and $S'$ of the same degree. Then, they have the same number of branch points and the same branch type, i.e., ${\rm type}(S)={\rm type}(S')$. In addition, if $\partial S=\partial S'$, then the isotopy can be chosen so that $h_u|_{D_1\times\partial D_2}=id$ for each $u\in[0,1]$. Thus, two surface braids are equivalent if and only if they are ambient isotopic by a fiber-preserving isotopy of $D_1\times D_2$ that fixes $D_1\times\partial D_2$ pointwise. \medskip

Kamada proved the following theorem (\cite{K1996, K-b}).

\begin{theorem} [\cite{K1996, K-b}]
Let $S$ and $S'$ be braided surfaces of degree $m$ with $n$ branch points whose boundary braids are the same $\delta\in B_m$. Let $\vec{b}$ and $\vec{b'}$ be braid systems of $S$ and $S'$, respectively, which are elements of $P^n_{\delta}(A_m)$. It holds that $S$ and $S'$ are equivalent if and only if  $\vec{b} \stackrel{{ \mathrm{Hur}}}{\backsim} \vec{b'}$.
\label{completely invariant of braided surface}
\end{theorem}
\medskip 

Take a braid $\delta\in B_m$ such that $P^n_{\delta}(A_m)\ne \emptyset$. By Theorems \ref{characterization of a braid system} and \ref{completely invariant of braided surface}, the Hurwitz equivalence class of braid systems in $P^n_{\delta}(A_m)$ is a complete invariant of braided surfaces of degree $m$ with $n$ branch points whose boundary braid is $\delta$. In particular, the Hurwitz equivalence class of braid systems in $P^n_{id}(A_m)$ is a complete invariant of surface braids of degree $m$ with $n$ branch points. Thus, we have the following.

\medskip 
\begin{corollary}Take a braid $\delta\in B_m$ such that $P^n_{\delta}(A_m)\ne \emptyset$. By restricting to $P^n_{\delta}(A_m)$, the Hurwitz invariants on $(B_m)^n$ obtained in Corollary \ref{cor-H-rank} are also invariants of braided surfaces of degree $m$ with $n$ branch points and boundary braid $\delta$. In particular, when $\delta=id$, they are invariants of surface braids of degree $m$ with $n$ branch points.
\label{cor-S-H-rank}
\end{corollary}
\medskip

\begin{remark} Through the natural homomorphism $f: B_m\to S_m$, the Hurwitz equivalence classes in ${(S_m)}^n$ are an invariant of braided surfaces of degree $m$ with $n$ branch points. Furthermore, the multiset of the cycle types of the $n$ permutations in an element $\vec{c}\in {(S_m)}^n$ coincides with the branch type of a braided surface which has a braid system $\vec{b}\in {(A_m)}^n$ such that $f^n(\vec{b})=\vec{c}$, and the branch type of a braided surface is also an invariant of braided surfaces. Prior to their role in the classification of braided surfaces, the Hurwitz equivalence classes in ${(S_m)}^n$, together with the branch types derived from them, have played fundamental roles in the classification of branched coverings of surfaces, as in \cite{Eis}. The invariants obtained in Corollary \ref{cor-H-rank} can distinguish braided surfaces with the same branch type. Indeed, the braid systems \(\vec b\) and \(\vec b'\) in Example \ref{ex-H} belong to \(P^4_{id}(A_4)\), and hence are braid systems of surface braids \(S\) and \(S'\) of degree \(4\) with 4 branch points. Their branch types coincide: ${\rm type}(S)={\rm type}(S')
=\{(4),(2,1,1),(2,1,1),(2,1,1)\}$. On the other hand, the multisets of characteristic polynomials for \(\vec b\) and \(\vec b'\) are different, as shown in Example \ref{ex-H}. Thus, \(S\) and \(S'\) are not equivalent, although they have the same branch type. We emphasize that our invariants effectively distinguish $S$ and $S'$ in this case.
\label{data of types of branch points}
\end{remark}
\medskip

\begin{remark}
Take a braided surface $S$ of degree $m$ with $n$ branch points and a Hurwitz generating system $(\gamma_1,\dots , \gamma_n)$ of $\pi_1(D_2\setminus \Sigma_S, y_0)$. 
We say $S$ is {\it simple} if ${\rm type}(S;y)=(2,1,\dots ,1)$ (i.e., $|p^{-1}(\{y\}) \cap S| = m-1$) for each $y\in \Sigma_S$. If $S$ is simple, then
$\rho_S(\gamma_i)\stackrel{\rm conj}{\backsim}\sigma_1$
or
$\rho_S(\gamma_i)\stackrel{\rm conj}{\backsim}\sigma_1^{-1}$
in $B_m$ for each $i\in\{1,\ldots,n\}$ (see \cite{K1996, K-b}). Let $SA_m$ be the set of braids $b\in B_m$ such that $b \stackrel{{ \mathrm{conj}}}{\backsim} \sigma_1$ or $b \stackrel{{ \mathrm{conj}}}{\backsim} \sigma_1^{-1}$. Note that $SA_m$ is a subset of $A_m$. Take a braid $\delta\in B_m$ such that $P^n_{\delta}(SA_m)\ne \emptyset$. The Hurwitz equivalence class of braid systems in $P^n_{\delta}(SA_m)$ is a complete invariant of simple braided surfaces of degree $m$ with $n$ branch points whose boundary braid is $\delta$. However, the Hurwitz invariants on $(B_m)^n$ obtained in Corollary \ref{cor-H-rank} provide no nontrivial invariants for simple braided surfaces; both $\sigma_1$ and $\sigma_1^{-1}$ have braid permutation $(1\ 2)\in S_m$, and $C(\sigma_1^2)$ (resp. $C(\sigma_1^{-2})$) is an $m\times m$ matrix whose nonzero entries are only the $(1,2)$- and $(2,1)$-entries, both equal to $1$ (resp. $-1$). By Proposition \ref{prop-conj-perm}, for each $b\in SA_m$, we have $C(b^2)\stackrel{\mathrm{perm}}{\backsim} C(\sigma_1^2)$ if $b\stackrel{\mathrm{conj}}{\backsim} \sigma_1$, and $C(b^2)\stackrel{\mathrm{perm}}{\backsim} C(\sigma_1^{-2})$ if $b\stackrel{\mathrm{conj}}{\backsim} \sigma_1^{-1}$. (For the non-simple case, these invariants are effective on $P^n_{\delta}(A_m)$ to distinguish certain non-simple braided surfaces, as illustrated in Example \ref{ex-H}.)
\label{simple braided surfaces}
\end{remark}
\medskip

\begin{remark}
Two braided surfaces of degree $m$ are {\it w-equivalent} if they are ambient isotopic by an isotopy $\{h_u\}_{u\in [0,1]}$ of $D_1\times D_2$ such that for each $u\in [0,1]$, $h_u$ is fiber-preserving and $h_1(Q_m\times \{y_0\})=Q_m\times \{y_0\}$ while $h_u$ preserves $p_2^{-1}(y_0)=D_1\times \{y_0\}$ setwise (not necessarily pointwise). Note that two equivalent braided surfaces are w-equivalent. Two w-equivalent braided surfaces of the same degree have the same number of branch points and the same branch type. Let $S$ and $S'$ be braided surfaces of degree $m$ with $n$ branch points and let $\vec{b}, \vec{b'}\in {(A_m)}^n$ be braid systems of $S$ and $S'$, respectively. Then, $S$ and $S'$ are w-equivalent if and only if $\vec{b} \stackrel{{ \mathrm{HC}}}{\backsim} \vec{b'}$, where $\vec{b} \stackrel{{ \mathrm{HC}}}{\backsim} \vec{b'}$ means that they are sent to each other by the Hurwitz action and the global conjugation, i.e., there exist elements $\beta\in B_n$ and $a\in B_m$ such that $(\vec{b}\cdot \beta) * a =\vec{b'}$. (The authors learned this assertion from Seiichi Kamada. According to him, the assertion can be proved by using the argument of Lemma 3 in \cite{K1996}.) By Proposition \ref{prop-rank}, the invariants of braided surfaces in Corollary \ref{cor-H-rank} up to equivalence are also invariants of them up to w-equivalence.
\label{w-equivalence}
\end{remark}
\medskip

\section{Application to surface links}
\label{section-surface-link}

An {\it oriented surface link} is a closed (orientable and) oriented surface embedded in ${\mathbb R}^4$ locally flatly. In this section, we simply call it a {\it surface link}. Two surface links $F$ and $F'$ are {\it equivalent} if there exists an orientation-preserving homeomorphism $h:{\mathbb R}^4\to {\mathbb R}^4$ such that $h(F)=F'$ and $h|_{F}:F\to F'$ is also an orientation-preserving homeomorphism. \medskip

Let $D_1$ and $D_2$ be 2-disks. We regard $D_1$ as a subset of a sphere $S^2$. Take a surface braid $S$ of degree $m$ in $D_1\times D_2$. Recall that $\partial S=Q_m\times \partial D_2^2$. The {\it closure} of $S$ is a closed oriented surface $\widehat{S}$ in $D_1\times S^2(\subset {\mathbb R}^4)$ such that $\widehat{S}\cap (D_1\times D_2)=S$ and $\widehat{S}\cap (D_1\times (S^2\setminus {\rm Int}(D_2)))=Q_m\times (S^2\setminus {\rm Int}(D_2))$. We regard the closed oriented surface $\widehat{S}$ as a surface link.\medskip

Kamada proved the following (\cite{K-b}).\medskip

\begin{theorem}[\cite{K-b}] The following hold.
\begin{itemize}
\item[(i)] ({\it Alexander's theorem in dimension four;})  Every surface link is equivalent to the closure of some surface braid. (In \cite{K-b}, an even stronger claim is proved as follows: every surface link is equivalent to the closure of some {\it simple surface braid} $S$, that is, a surface braid such that ${\rm type}(S;y)=(2,1,\dots, 1)$ for every branch point $y\in \Sigma_S$.)
\item[(ii)] ({\it Markov’s theorem in dimension four;}\footnote{Application of the crossing matrix regarding the classical Markov's theorem is discussed in \cite{AS}.}) The closures of two surface braids $S$ and $S'$ are equivalent as surface links if and only if $S$ and $S'$ are related by the deformations ``braid ambient isotopic'', ``conjugations'', ``stabilizations'', and their inverse operations. 
\end{itemize}
\label{A&Mtheorem}
\end{theorem}



\noindent In terms of the braid systems, we have the following theorem. 

\medskip 
\begin{theorem}[\cite{K-b}]
Let $m, m', n, n'\in {\mathbb N}$. 
Let $\vec{b}\in P^n_{id}(A_m)$ and $\vec{b'}\in P^{n'}_{id}(A_{m'})$. 
Let $S$ and $S'$ be surface braids which have braid systems $\vec{b}$ and $\vec{b'}$, respectively. 
If their closures $\widehat{S}$ and $\widehat{S'}$ are equivalent as surface links, then $\vec{b}$ and $\vec{b'}$ are related by a finite sequence of the following transformations (I)--(I\hspace{-0.5pt}V'): For $k, l\in {\mathbb N}$, 
\begin{itemize}
\item[(I)] applying the Hurwitz action by a braid in $B_k$ on $P^l_{id}(A_k)$, 
\item[(I\hspace{-0.5pt}I)] applying the global conjugation by a braid in $B_k$ on $P^l_{id}(A_k)$, 
\item[(I\hspace{-0.5pt}I\hspace{-0.5pt}I)] applying a stabilization $P^l_{id}(A_k)\to P^{l+2}_{id}(A_{k+1})$,
\item[(I\hspace{-0.5pt}I\hspace{-0.5pt}I')] applying a destabilization $P^{l+2}_{id}(A_{k+1})\to P^{l}_{id}(A_{k})$,
\item[(I\hspace{-0.5pt}V)] applying an Euler fission $P^{l}_{id}(A_k)\to P^{l+q}_{id}(A_k)$ for $q>0$, that transforms $(b_1,\dots ,b_l)\in P^{l}_{id}(A_k)$ to $(b_1',\dots ,b_{l+q}')\in P^{l+q}_{id}(A_k)$ with $b_i=b_{i'}$ for $1\leq i\leq l-1$, $b_l=b_l'b_{l+1}' \cdots b_{l+q}'$ in $B_k$ and $\tau (b_l)= \tau (b_l') + \tau (b_{l+1}') + \dots + \tau (b_{l+q}')$, where $\tau (b)$ denotes the difference $k$ minus the component number of the closure of $b$ for $b\in B_k$, and 
\item[(I\hspace{-0.5pt}V')] applying an Euler fusion $P^{l+q}_{id}(A_k)\to P^{l}_{id}(A_k)$ for $q>0$, that is the inverse of an Euler fission. 
\end{itemize}
\label{thm-1234}
\end{theorem}
\medskip 

\noindent The transformations (I), (I\hspace{-0.5pt}I), (I\hspace{-0.5pt}I\hspace{-0.5pt}I) and (I\hspace{-0.5pt}V) in Theorem \ref{thm-1234} correspond to ``equivalence'', ``conjugations'', ``stabilizations'' and ``braid ambient isotopic'' on surface braids, respectively, in the property (ii) in Theorem \ref{A&Mtheorem}. \medskip

\noindent In the following example, the essential eigenvalue set $E(\vec{b})$ detects the necessity of the Euler fission (I\hspace{-0.5pt}V) or fusion (I\hspace{-0.5pt}V')\footnote{For the case of Example \ref{ex-E}, we can also detect the necessity of the Euler fusion or fission between $\vec{b}$ and $\vec{c}$ by comparing the degrees and the lengths; the value of $2k-l$ is unchanged by the operations (I), (I\hspace{-0.5pt}I), (I\hspace{-0.5pt}I\hspace{-0.5pt}I), and (I\hspace{-0.5pt}I\hspace{-0.5pt}I') for the degree $k$ and the length $l$ of $\vec{b} \in P^l_{id}(A_k)$. 
At the moment, pairs of braid systems $\vec{b}$ and $\vec{c}$ with the same value of $2k-l$ and $E(\vec{b}) \neq E(\vec{c})$ have not been found.}. 

\medskip 
\begin{example}
For a braid system $\vec{b}=(\sigma_1 \sigma_2^{-1} \sigma_3, \ \sigma_3^{-1}, \ \sigma_2, \ \sigma_1^{-1}) \in P^4_{id}(A_4)$, we obtain the set $E(\vec{b})= \{ -3 \}$ from the characteristic polynomial $P(\vec{b})=x^6 (x+1)^3 (x-1)^6 (x+3)$. 
For the braid system $\vec{c}=(\sigma_1 \sigma_2^{-1} \sigma_3, \ \sigma_3^{-1}\sigma_2 \sigma_1^{-1}) \in P^2_{id}(A_4)$, we obtain the set $E(\vec{c})= \{ \pm 3 \}$ from the characteristic polynomial $P(\vec{c})=(x+1)^3 (x-1)^3 (x+3)(x-3)$. 
Note that $\vec{c}$ is obtained from $\vec{b}$ by applying Euler fusions (I\hspace{-0.5pt}V') twice. 
Moreover, the difference between $E(\vec{b})= \{ -3 \}$ and $E(\vec{c})= \{ \pm 3 \}$ implies that any sequence of (I)--(I\hspace{-0.5pt}V') between $\vec{b}$ and $\vec{c}$ must contain an Euler fusion or fission by Theorem \ref{thm-main-2}. 
\label{ex-E}
\end{example}

\section*{Acknowledgment}
The authors are deeply grateful to Seiichi Kamada for valuable comments and suggestions, in particular regarding applications to braided surfaces, surface braids, and surface links. 
This work was partially supported by the JSPS KAKENHI Grant Number JP21K03263.

\end{document}